\documentclass{ws-procs975x65-ppn}

\newcommand{\R}{{\mathbb R}}
\renewcommand{\S}{{\mathbb S}^{d-1}}
\newcommand{\N}{{\mathbb N}}

\newcommand{\ird}[1]{\int_{\R^d}{#1}\;dx}
\newcommand{\nrm}[2]{\|{#1}\|_{\L^{#2}(\R^d)}}
\newcommand{\nrmcnd}[2]{\|{#1}\|_{\L^{#2}(\mathcal C)}}
\newcommand{\icnd}[1]{\int_{\mathcal C}{#1}\;dy}
\newcommand{\be}[1]{\begin{equation}\label{#1}}
\newcommand{\ee}{\end{equation}}
\renewcommand{\(}{\left(}
\renewcommand{\)}{\right)}
\renewcommand{\S}{{\mathbb S^{d-1}}}
\newcommand{\C}[1]{\mathsf C_{\rm #1}}
\renewcommand{\H}{\mathrm H}
\renewcommand{\L}{\mathrm L}

\newcommand{\CKN}{{(CKN)}}
\newcommand{\WLH}{{(WLH)}}

\begin{document}
\thispagestyle{empty}

\title{Extremal functions in some interpolation inequalities:\\
Symmetry, symmetry breaking and estimates of the best constants}

\author{J. Dolbeault$^*$ and M.J. Esteban$^{**}$}

\address{Ceremade (UMR CNRS no.~7534), Universit\'e Paris-Dauphine,\\ Place de Lattre de Tassigny, F-75775 Paris C\'edex 16, France\\$^*$E-mail: dolbeaul@ceremade.dauphine.fr, $^{**}$E-mail: esteban@ceremade.dauphine.fr\\
http://www.ceremade.dauphine.fr/$\sim$dolbeaul/,
http://www.ceremade.dauphine.fr/$\sim$esteban/}

\begin{abstract}
This contribution is devoted to a review of some recent results on existence, symmetry and symmetry breaking of optimal functions for Caffarelli-Kohn-Nirenberg \CKN\ and weighted logarithmic Hardy \WLH\ inequalities. These results have been obtained in a series of papers \cite{DET,0902,DDFT,DE2010,1007} in collaboration with M.~del Pino, S.~Filippas, M.~Loss, G.~Tarantello and A.~Tertikas and are presented from a new viewpoint.
\end{abstract}

\keywords{Caffarelli-Kohn-Nirenberg inequality; Gagliardo-Nirenberg inequality; logarithmic Hardy inequality; logarithmic Sobolev inequality; extremal functions; radial symmetry; symmetry breaking; Emden-Fowler transformation; linearization; existence; compactness; optimal constants}

\bodymatter

\section{Two families of interpolation inequalities}\label{Sec:Interpolation}

Let $d\in\N^*$, $\theta\in[0,1]$, consider the set $\mathcal D$ of all smooth functions which are compactly supported in $\R^d\setminus\{0\}$ and define $\vartheta(d,p):=d\,\frac{p-2}{2\,p}$, $a_c:=\frac{d-2}2$, $\Lambda(a):=(a-a_c)^2$ and $p(a,b):=\frac{2\,d}{d-2+2\,(b-a)}$. We shall also set $2^*:=\frac{2\,d}{d-2}$ if $d\ge 3$ and $2^*:=\infty$ if $d=1$ or~Ê$2$. For any $a<a_c$, we consider the two families of interpolation inequalities:
\begin{description}
\item[(CKN)] \emph{Caffarelli-Kohn-Nirenberg inequalities} \cite{Caffarelli-Kohn-Nirenberg-84,DDFT,DE2010} -- Let $b\in(a+1/2,a+1]$ and $\theta\in(1/2,1]$ if $d=1$, $b\in(a,a+1]$ if $d=2$ and $b\in[a,a+1]$ if $d\ge3$. Assume that $p=p(a,b)$, and $\theta\in[\vartheta(d,p),1]$ if $d\ge2$. There exists a finite positive constant $\C{CKN}(\theta,p,a)$ such that, for any $u\in\mathcal D$,
\[
\nrm{|x|^{-b}\,u}p^2\le\C{CKN}(\theta,p,a)\,\nrm{|x|^{-a}\,\nabla u}2^{2\,\theta}\,\nrm{|x|^{-(a+1)}\,u}2^{2\,(1-\theta)}\,.
\]
\item[(WLH)] \emph{Weighted logarithmic Hardy inequalities} \cite{DDFT,DE2010} --
Let $\gamma \ge d/4$ and $\gamma>1/2$ if $d=2$. There exists a positive constant $\C{WLH}(\gamma,a)$ such that, for any $u\in\mathcal D$, normalized by $\nrm{|x|^{-(a+1)}\,u}2=1$,
\[
\ird{\frac{|u|^2\,\log\(|x|^{d-2-2\,a}\,|u|^2\)}{|x|^{2\,(a+1)}}}\leq2\,\gamma\,\log\!\left[\C{WLH}(\gamma,a)\,\nrm{\,|x|^{-a}\,\nabla u}2^2\right]\,.
\]
\end{description}
\WLH\ appears as a limiting case \cite{DDFT,DE2010} of \CKN\ with $\theta=\gamma\,(p-2)$ as $p\to 2_+$. By a standard completion argument, these inequalities can be extended to the set $\mathcal D_a^{1,2}(\R^d):=\{u\in\L^1_{\rm loc}(\R^d)\,:\,|x|^{-a}\,\nabla u\in\L^2(\R^d)\;\mbox{\emph{and}}\;|x|^{-(a+1)}\, u\in\L^2(\R^d)\}$. We shall assume that all constants in the inequalities are taken with their optimal values. For brevity, we shall call \emph{extremals} the functions which realize equality in \CKN\ or in \WLH.

Let $\C{CKN}^*(\theta,p,a)$ and $\C{WLH}^*(\gamma,a)$ denote the optimal constants when admissible functions are restricted to the radial ones. \emph{Radial extremals} are explicit and the values of the constants, $\C{CKN}^*(\theta,p,a)$ and $\C{WLH}^*(\gamma,a)$, are known \cite{DDFT}. Moreover,~we~have
\be{Ineq:CompRad}\begin{array}{l}
\C{CKN}(\theta,p,a)\ge\C{CKN}^*(\theta,p,a)=\C{CKN}^*(\theta,p,a_c-1)\,\Lambda(a)^{\frac{p-2}{2p}-\theta}\;,\\[6pt]
\C{WLH}(\gamma,a)\ge\C{WLH}^*(\gamma,a)=\C{WLH}^*(\gamma,a_c-1)\,\Lambda(a)^{-1+\frac 1{4\,\gamma}}\;.
\end{array}\ee
Radial symmetry for the extremals of \CKN\ and \WLH\ implies that $\C{CKN}(\theta,p,a)=\C{CKN}^*(\theta,p,a)$ and $\C{WLH}(\gamma,a)=\C{WLH}^*(\gamma,a)$, while  \emph{symmetry breaking} only means that inequalities in \eqref{Ineq:CompRad} are strict.

\section{Existence of extremals}\label{Sec:Extremals}

\begin{theorem}\label{Thm:Existence} Equality \cite{DE2010} in \CKN\ is attained for any $p\in(2,2^*)$ and $\theta\in(\vartheta(p,d),1)$ or $\theta=\vartheta(p,d)$ and $a\in(a_\star^{\scriptscriptstyle\rm CKN},a_c)$, for some $a_\star^{\scriptscriptstyle\rm CKN}<a_c$. It is not attained if $p=2$, or $a<0$, $p=2^*$, $\theta=1$ and $d\ge 3$, or $d=1$ and $\theta=\vartheta(p,1)$.

Equality \cite{DE2010} in \WLH\ is attained if $\gamma\ge 1/4$ and $d=1$, or $\gamma>1/2$ if $d=2$, or for $d\geq 3$ and either $\gamma>d/4\,$ or $\,\gamma=d/4$ and $a\in(a_\star^{\scriptscriptstyle\rm WLH},a_c)$, where $a_\star^{\scriptscriptstyle\rm WLH}:=a_c-\sqrt{\Lambda_\star^{\scriptscriptstyle\rm WLH}}$ and $\Lambda_\star^{\scriptscriptstyle\rm WLH}:=(d-1)\,e\,(2^{d+1}\,\pi)^{-1/(d-1)}\,\Gamma(d/2)^{2/(d-1)}$.\end{theorem}

Let us give some hints on how to prove such a result. Consider first Gross' logarithmic Sobolev inequality in Weissler's form \cite{MR479373}
\[
\ird{|u|^2\,\log |u|^2}\le\frac d2\,\log\(\C{LS}\,\nrm{\nabla u}2^2\)\quad\forall\; u\in\H^1(\R^d)\;\mbox{s.t.}\;\nrm u2=1\;.
\]
The function $u(x)=(2\,\pi)^{-d/4}\,\exp(-|x|^2/4)$ is an extremal for such an inequality. By taking $u_n(x):=u(x+n\,\mathsf e)$ for some $\mathsf e\in\S$ and any $n\in\N$ as test functions for \WLH, and letting $n\to+\infty$, we find that $\C{LS}\le\C{WLH}(d/4,a)$. If equality holds, this is a mechanism of loss of compactness for minimizing sequences. On the opposite, if $\C{LS}<\C{WLH}(d/4,a)$, which is the case if $a\in(a_\star^{\scriptscriptstyle\rm WLH},a_c)$ where $a_\star^{\scriptscriptstyle\rm WLH}=a$ is given by the condition $\C{LS}=\C{WLH}^*(d/4,a)$, we can establish a compactness result which proves that equality is attained in \WLH\ in the critical case $\gamma=d/4$.

A similar analysis  for \CKN\ shows that $\C{GN}(p)\le\C{CKN}(\theta,p,a)$ in the critical case $\theta=\vartheta(p,d)$, where $\C{GN}(p)$ is the optimal constant in the Gagliardo-Nirenberg-Sobolev interpolation inequalities
\[
\nrm up^2\le\C{GN}(p)\,\nrm{\nabla u}2^{2\,\vartheta(p,d)}\,\nrm u2^{2\,(1-\vartheta(p,d))}\quad\forall\; u\in\H^1(\R^d)
\]
and $p\in(2,2^*)$ if $d=2$ or $p\in(2,2^*]$ if $d\ge3$. However, extremals are not known explicitly in such inequalities if $d\ge2$, so we cannot get an explicit interval of existence in terms of $a$, even if we also know that compactness of minimizing sequences for \CKN\ holds when  $\C{GN}(p)<\C{CKN}(\vartheta(p,d),p,a)$. This is the case if $a>a_\star^{\scriptscriptstyle\rm CKN}$ where $a=a_\star^{\scriptscriptstyle\rm CKN}$ is defined by the condition $\C{GN}(p)=\C{CKN}^*(\vartheta(p,d),p,a)$.

It is very convenient to reformulate \CKN\ and \WLH\ inequalities in cylindrical variables \cite{Catrina-Wang-01}. By means of the Emden-Fowler transformation
\[
s=\log|x|\in\R\;,\quad\omega=x/|x|\in\S\,,\quad y=(s,\omega)\;,\quad v(y)=|x|^{a_c-a}\,u(x)\;,
\]
\CKN\ for $u$ is equivalent to a Gagliardo-Nirenberg-Sobolev inequality on the cylinder $\mathcal C:=\R\times\S$ for $v$, namely
\[
\nrmcnd vp^2\leq\C{CKN}(\theta,p,a)\(\;\nrmcnd{\nabla v}2^2+\Lambda\,\nrmcnd v2^2\)^\theta\,\nrmcnd v2^{2\,(1-\theta)}\quad\forall\;v\in\H^1(\mathcal C)
\]
with $\Lambda=\Lambda(a)$. Similarly, with $w(y)=|x|^{a_c-a}\,u(x)$, \WLH\ is equivalent to
\[\label{Ineq:GLogHardy-w}
\icnd{|w|^2\,\log |w|^2}\leq 2\,\gamma\,\log\left[\C{WLH}(\gamma,a)\left(\nrmcnd{\nabla w}2^2+\Lambda\right)\right]
\]
for any $w\in\H^1(\mathcal C)$ such that $\nrmcnd w2=1$. Notice that radial symmetry for $u$ means that $v$ and $w$ depend only on $s$.

Consider a sequence $(v_n)_n$ of functions in $\H^1(\mathcal C)$, which minimizes the functional
\[
\mathcal E_{\theta,\Lambda}^p[v]:=\(\nrmcnd{\nabla v}2^2\!+\Lambda\,\nrmcnd v2^2\)^\theta\nrmcnd v2^{2\,(1-\theta)}
\]
under the constraint $\nrmcnd{v_n}p=1$ for any $n\in\N$. As quickly explained below, if bounded, such a sequence is relatively compact and converges up to translations and the extraction of a subsequence towards a minimizer of $\mathcal E_{\theta,\Lambda}^p$.

Assume that $d\ge3$, let $t:=\nrmcnd{\nabla v}2^2/\nrmcnd v2^2$ and $\Lambda=\Lambda(a)$. If $v$ is a minimizer of $\mathcal E_{\theta,\Lambda}^p[v]$ such that $\nrmcnd vp=1$, then we have
\[
(t+\Lambda)^\theta=\mathcal E_{\theta,\Lambda}^p[v]\,\frac{\nrmcnd vp^2}{\nrmcnd v2^2}=\frac{\nrmcnd vp^2}{\C{CKN}(\theta,p,a)\,\nrmcnd v2^2}\le\frac{\mathsf S_d^{\vartheta(d,p)}}{\C{CKN}(\theta,p,a)}\,\(t\!+a_c^2\)^{\vartheta(d,p)}
\]
where $\mathsf S_d=\C{CKN}(1,2^*,0)$ is the optimal Sobolev constant, while we know from \eqref{Ineq:CompRad} that $\lim_{a\to a_c}\C{CKN}(\theta,p,a)=\infty$ if $d\ge 2$. This provides a bound on $t$ if $\theta>\vartheta(p,d)$. An estimate can be obtained also for $v_n$, for $n$ large enough, and standard tools of the concentration-compactness method allow to conclude that $(v_n)_n$ converges towards an extremal. A similar approach holds for \CKN\ if $d=2$, or for \WLH.

The above variational approach also provides an existence result of extremals for \CKN\ in the critical case $\theta=\vartheta(p,d)$, if $a\in(a_1,a_c)$ where $a_1:=a_c-\sqrt{\Lambda_1}$ and $\Lambda_1=\min\{(\C{CKN}^*(\theta,p,a_c-1)^{1/\theta}/\,\mathsf S_d)^{d/(d-1)},(a_c^2\,\C{CKN}^*(\theta,p,a_c-1)^{1/\theta}/\,\mathsf S_d)^d$.

\medskip If symmetry is known, then there are (radially symmetric) extremals \cite{DDFT}. Anticipating on the results of the next section, we can state the following result which arises as a consequence of Schwarz' symmetrization method (see Theorem~\ref{Thm:Symmetry}, below).
\begin{proposition}\label{Prop} Let $d\ge3$. Then \CKN\ with $\theta=\vartheta(p,d)$ admits a radial extremal if \cite{1007} $a\in[a_0,a_c)$ where $a_0:=a_c-\sqrt{\Lambda_0}$ and $\Lambda=\Lambda_0$ is defined by the condition $\Lambda^{(d-1)/d}=\vartheta(p,d)\,\C{CKN}^*(\theta,p,a_c-1)^{1/\vartheta(d,p)}/\,\mathsf S_d$. \end{proposition}
A similar estimate also holds if $\theta>\vartheta(d,p)$, with less explicit computations.\cite{1007}

\section{Symmetry and symmetry breaking}\label{Sec:Symmetry}

Define
\[\begin{array}{l}
\underline a(\theta,p):= a_c-\frac{2\,\sqrt{d-1}}{p+2}\,\sqrt{\frac{2\,p\,\theta}{p-2}-1}\;,\quad \tilde a(\gamma):=a_c-\frac 12\sqrt{(d-1)(4\,\gamma-1)}\;,\\
\Lambda_{\rm SB}(\gamma):=\frac 18\,(4\,\gamma-1)\,e\,\big(\tfrac{\pi^{4\,\gamma-d-1}}{16}\big)^\frac 1{4\,\gamma-1}\,\big(\tfrac d\gamma\big)^\frac{4\,\gamma}{4\,\gamma-1}\,\Gamma\(\tfrac d2\)^\frac 2{4\,\gamma-1}\,.
\end{array}\]
\begin{theorem}\label{Thm:SymmetryBreaking} Let $d\ge 2$ and $p\in(2,2^*)$. Symmetry breaking holds in \CKN\ if either \cite{DDFT,1007} $a<\underline a(\theta,p)$ and $\theta\in[\vartheta(p,d),1]$, or \cite{1007} $a<a_\star^{\scriptscriptstyle\rm CKN}$ and $\theta=\vartheta(p,d)$.

Assume that $\gamma>1/2$ if $d=2$ and $\gamma\ge d/4$ if $d\ge3$. Symmetry breaking holds in \WLH\ if \cite{DDFT,1007} $a<\max\{\tilde a(\gamma), a_c-\sqrt{\Lambda_{\rm SB}(\gamma)}\}$.\end{theorem}
When $\gamma=d/4$, $d\ge 3$, we observe that $\Lambda_\star^{\scriptscriptstyle\rm WLH}=\Lambda_{\rm SB}(d/4)<\Lambda(\tilde a(d/4))$ with the notations of Theorem~\ref {Thm:Existence} and there is symmetry breaking if $a\in(-\infty,a_\star^{\scriptscriptstyle\rm WLH})$, in the sense that $\C{WLH}(d/4,a)>\C{WLH}^*(d/4,a)$, although we do not know if extremals for \WLH\ exist when $\gamma=d/4$.

Results of symmetry breaking for \CKN\ with $a<\underline a(\theta,p)$ have been established first \cite{Catrina-Wang-01,Felli-Schneider-03,DET} when $\theta=1$ and later\cite{DDFT} extended to $\theta<1$. The main idea in case of \CKN\ is consider the quadratic form associated to the second variation of $\mathcal E_{\theta,\Lambda}^p$ around a minimizer among functions depending on $s$ only and observe that the linear operator $\mathcal L_{\theta,\Lambda}^p$ associated to the quadratic form has a negative eigenvalue if $a<\underline a$. Results \cite{DDFT} for \WLH, $a<\tilde a(\gamma)$, are based on the same method.

For any $a<a_\star^{\scriptscriptstyle\rm CKN}$, we have $\C{CKN}^*(\vartheta(p,d),p,a)<\C{GN}(p)\le\C{CKN}(\vartheta(p,d),p,a)$, which proves symmetry breaking. Using well-chosen test functions, it has been proved \cite{1007} that $\underline a(\vartheta(p,d),p)<a_\star^{\scriptscriptstyle\rm CKN}$ for $p-2>0$, small enough, thus also proving symmetry breaking for $a-\underline a(\vartheta(p,d),p)>0$, small, and $\theta-\vartheta(p,d)>0$, small.

\begin{theorem}\label{Thm:Symmetry} For all $d\geq 2$, there exists \cite{0902,1007} a continuous function $a^*$ defined on the set $\{(\theta,p)\in(0,1]\times(2,2^*)\,:\,\theta>\vartheta(p,d)\}$ such that $\lim_{p\to 2_+}a^*(\theta,p)=-\infty$ with the property that \CKN\ has only radially symmetric extremals if $(a,p)\in(a^*(\theta,p),a_c)\times(2,2^*)$, and none of the extremals is radially symmetric if $(a,p)\in(-\infty,a^*(\theta,p))\times(2,2^*)$.

Similarly, for all $d\geq 2$, there exists \cite{1007} a continuous function $a^{**}:(d/4,\infty)\to(-\infty,a_c)$ such that, for any $\gamma>d/4$ and $a\in [a^{**}(\gamma), a_c)$, there is a radially symmetric extremal for \WLH, while for $a<a^{**}(\gamma)$ no extremal is radially symmetric.
\end{theorem}
Schwarz' symmetrization allows to characterize \cite{1007} a subdomain of $(0,a_c)\times(0,1)\ni(a,\theta)$ in which symmetry holds for extremals of \CKN, when $d\ge 3$. If $\theta=\vartheta(p,d)$ and $p>2$, there are radially symmetric extremals\cite{1007} if $a\in[a_0,a_c)$ where $a_0$ is given in Propositions~\ref{Prop}.

Symmetry also holds if $a-a_c$ is small enough, for \CKN\ as well as for \WLH, or when $p\to 2_+$ in \CKN, for any $d\ge 2$, as a consequence of the existence of the spectral gap of $\mathcal L_{\theta,\Lambda}^p$ when $a>\underline a(\theta,p)$.

For given $\theta$ and $p$, there is\cite{0902,1007} a unique $a^*\in(-\infty,a_c)$ for which there is symmetry breaking in $(-\infty,a^*)$ and for which all extremals are radially symmetric when $a\in(a^*,a_c)$. This follows from the observation that, if $v_\sigma(s,\omega):=v(\sigma\,s,\omega)$ for $\sigma>0$, then $(\mathcal E_{\theta,\sigma^2\Lambda}^p[v_\sigma])^{1/\theta}-\sigma^{(2\,\theta-1+2/p)/\theta^2}\,(\mathcal E_{\theta,\Lambda}^p[v])^{1/\theta}$ is equal to $0$ if $v$ depends only on~$s$, while it has the sign of $\sigma-1$ otherwise.

From Theorem~\ref{Thm:SymmetryBreaking}, we can infer that radial and non-radial extremals for \CKN\ with $\theta>\vartheta(p,d)$ coexist on the threshold, in some cases.

\medskip Numerical results illustrating our results on existence and on symmetry / symmetry breaking have been collected in Fig.~1 below in the critical case for \CKN.

\def\figsubcap#1{\par\noindent\centering\footnotesize(#1)}
\begin{figure}[!ht]
\begin{center}
\parbox{2.1in}{\epsfig{figure= 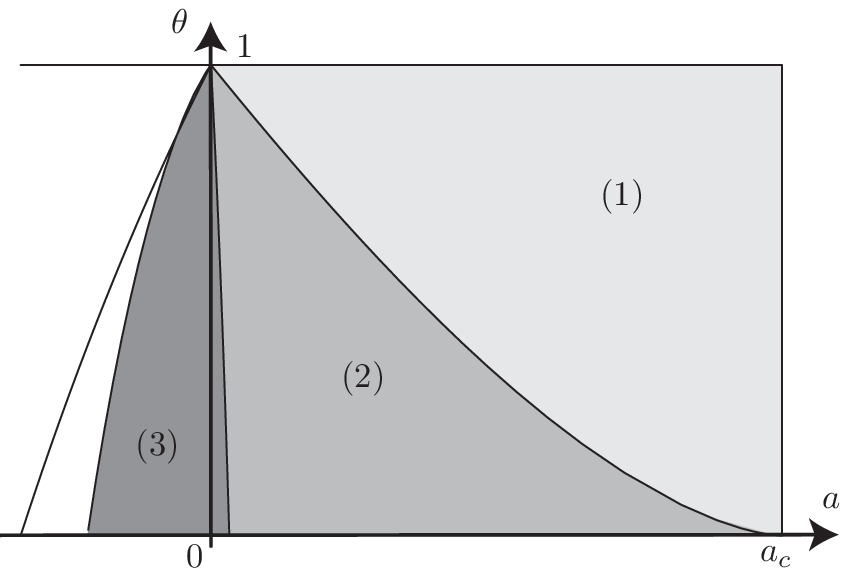,width=2in}\figsubcap{a}}
\hspace*{12pt}
\parbox{2.1in}{\epsfig{figure= 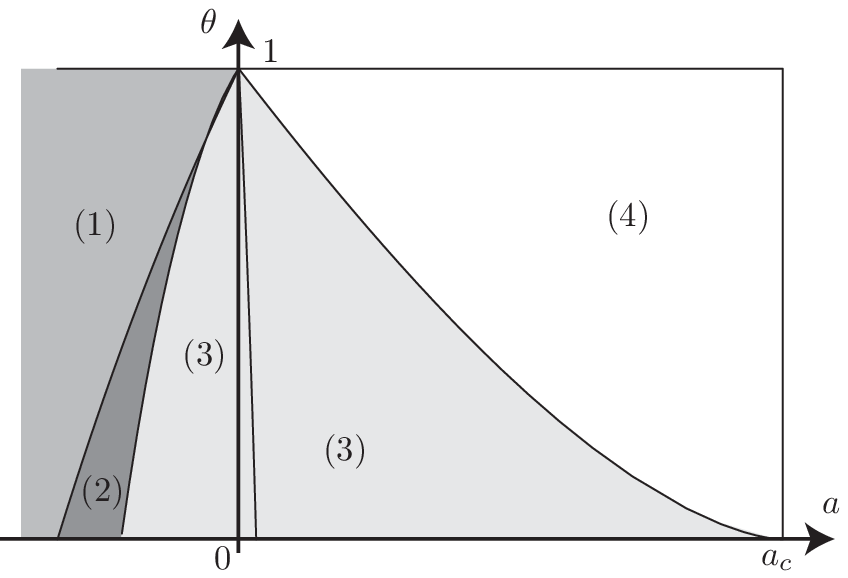,width=2in}\figsubcap{b}}
\caption{Critical case for \CKN: $\theta=\vartheta(p,d)$. Here we assume that $d=5$.\newline
(a) The zones in which existence is known are (1) in which $a\ge a_0$, because extremals are achieved among radial functions, (2) using the \emph{a priori} estimates: $a>a_1$, and (3) by comparison with the Gagliardo-Nirenberg inequality: $a>a_\star^{\scriptscriptstyle\rm CKN}$.\newline
(b) The zone of symmetry breaking contains (1) by linearization around radial extremals: $a<\underline a(\theta,p)$, and (2) by comparison with the Gagliardo-Nirenberg inequality: $a< a_\star^{\scriptscriptstyle\rm CKN}$; in (3) it is not known whether symmetry holds or if there is symmetry breaking, while in (4) symmetry holds by Schwarz' symmetrization: $a_0\le a<a_c$.\newline
Numerically, we observe that $\underline a$ and $a_\star^{\scriptscriptstyle\rm CKN}$ intersect for some $\theta\approx 0.85.$
}%
\label{fig1.2}
\end{center}
\end{figure}

\vspace*{-12pt}\noindent{\scriptsize\emph{Acknowledgements}.
The authors have been supported by the ANR projects CBDif-Fr and EVOL.\\[-4pt]\copyright\,2010 by the authors. This paper may be reproduced, in its entirety, for non-commercial purposes.}\vspace*{-6pt}

\end{document}